\documentclass[12pt,a4paper]{amsart}
\usepackage{amsfonts}
\usepackage{amsthm}
\usepackage{amsmath}
\usepackage{amscd}
\usepackage[utf8]{inputenc}
\usepackage{t1enc}
\usepackage[mathscr]{eucal}
\usepackage{indentfirst}
\usepackage{graphicx}
\usepackage{graphics}
\usepackage{pict2e}
\usepackage{epic}
\numberwithin{equation}{section}
\usepackage[margin=2.9cm]{geometry}
\usepackage{epstopdf}

\usepackage{MnSymbol}%afficher blacksquare

\theoremstyle{plain}
\newtheorem{Th}{Theorem}[section]
\newtheorem{Lemma}[Th]{Lemma}
\newtheorem{Cor}[Th]{Corollary}

 \theoremstyle{definition}
\newtheorem{Def}[Th]{Definition}

\newtheorem{Rem}[Th]{Remark}
\newtheorem{?}[Th]{Problem}

\begin{document}

\title{BEURLING'S  THEOREM FOR THE TWO-SIDED QUATERNION FOURIER TRANSFORM}

\author[S. Fahlaoui]{Said Fahlaoui}

\address{Said Fahlaoui \\ Department of Mathematics and Computer Sciences, Faculty of Sciences, Equipe d'Analyse Harmonique et Probabilit\'es, University Moulay Ismail, BP 11201 Zitoune, Meknes, Morocco}
\email{saidfahlaoui@gmail.com}

\author[Y. Elhaoui]{Youssef El Haoui}

\address{Youssef El Haoui\\ Department of Mathematics and Computer Sciences, Faculty of Sciences, Equipe d'Analyse Harmonique et Probabilit\'es, University Moulay Ismail, BP 11201 Zitoune, Meknes, Morocco}
\email{youssefelhaoui@gmail.com}

 \subjclass[2010]{}


\begin{abstract} 
 The  two-sided quaternion Fourier transform satisfies some uncertainty principles similar to the Euclidean Fourier transform.  A generalization  of  Beurling’s theorem, Hardy, Cowling-Price and Gelfand-Shilov theorems, is obtained for the two-sided quaternion Fourier transform.

\end{abstract}

\maketitle
\section{Introduction}
In harmonic analysis, the uncertainty principle states that a non zero function and its Fourier transform cannot both be sharply localized. This fact is expressed by several versions which  were proved by Hardy, Cowling-Price and Gelfand-Shilov ..[7,13]. A more general version of uncertainty principle, which is called Beurling’s theorem which is given by A. Beurling and proved by H\"ormander [11] and generalized by Bonami et al [2], asserts that

\begin{Th}   
Let  $f \in L^2({\mathbb R}^n)$ and  $d\ge 0\ $ such that

$\int_{{\mathbb R}^n}{\int_{\mathbb R^n}{\frac{\left|f(x)\right|\left|\hat{f}(y)\right|}{{(1+\left|x\right|+\left|y\right|)}^d}}}\ e^{2\pi \left|x\right|\left|y\right|}dxdy<\infty $   where $\hat{f}\left(y\right)=\int_{{\mathbb R}^n}{e^{-2\pi i<x,y>}}\ f(x) dx,$

Then $f(x)= P(x)\ e^{-<Ax,x>}, $

where $A$ is a real positive definite symmetric matrix and P is a polynomial of degree $<\frac{(d-n)}{2}$. In particular , $f$ is identically 0 when $d \le n.$

 \end{Th}

\vspace*{1 cm}

 Our paper is organized as follows. In section 2, we review basic notions and notations related to the quaternion algebra. In section 3, we recall the notion and

some results for the two-sided quaternionic Fourier transform useful in the sequel. In section 4, we prove the Beurling's theorem for the two-sided

quaternionic Fourier transform. Section 5 contains other uncertainty principles for  the two-sided QFT : Hardy and  Gelfand-Shilov.

Note that we will often use the shorthand  $x:= (x_1,x_2)$,$\ y:= (y_1,y_2)$, also the letter $ C$ indicates a positive constant that is not necessarily the same in each occurrence.

\section{The algebra of quaternions }
The quaternion algebra over $\mathbb{R}$, denoted by $\mathbb{H}$, is an associative noncommutative four-dimensional algebra, it was invented by  W. R. Hamilton in 1843.\\
$\mathbb{H}=\{q=q_0+iq_1+jq_2+kq_3;\ q_0,q_1,\ q_2,q_3 \in \mathbb{R}$\}
where ${  i},{  \ }{  j},{  \ }{  k}$ satisfy Hamilton's multiplication rules

$ij=-ji=k~ ; jk=-kj=i~ ; ki=-ik=j~ ; i^2 = j^2 = k^2=-1.$

Quaternions are isomorphic to the Clifford algebra ${Cl}_{(0,2)}$ of ${\mathbb R}^{(0,2)}$:

\hspace*{4 cm}$\mathbb{H} \cong {Cl}_{(0,2)}.$          \hfill(2.1)

The scalar part of  a quaternion  $q \in  \mathbb{H}$\ is $q_0$ denoted by $Sc(q)$, the non scalar part(or pure quaternion) of $q$\ is $iq_1+jq_2+kq_3$ denoted by $Vec(q)$.

We define the conjugation of $q \in \mathbb{H}$\ by ~:

$\overline{q}$=$q_0-iq_1-jq_2-kq_3$. \\ The quaternion conjugation is a linear anti-involution 

\[\overline{qp}= \overline{p} \ \overline{q} ,\ \overline{p+q}= \overline{p}+\overline{q},\ \overline{\overline{p}}=p.\] 

The modulus of a quaternion  q is defined by:$\ $

\hspace*{4 cm}${|q|}_Q=\sqrt{q\overline{q}}=\sqrt{{q_0}^2{  +}\ {q_1}^2{  +}{q_2}^2+{q_3}^2}.$ \hfill(2.2)

In particular, when $q=q_0$ is a real number, the module ${|q|}_Q$ reduces to the ordinary Euclidean module $\left|q\right|=\sqrt{{q_0}^2}$.

we have~:

\hspace*{4cm}${|pq|}_Q={|p|}_Q{|q|}_Q.$ \hfill(2.3)

It is easy to verify that  0$\ne q \in \mathbb{H}$ implies~:

\[q^{-1}=\frac{\overline{q}}{{{|q|}_Q}^2}.\] 

Any quaternion  $q$ can be written as $q$=\ ${{|q|}_Qe}^{\mu \theta }$  where $e^{\mu \theta }$ is understood in accordance with Euler's formula $e^{\mu \theta }={\cos  \left(\theta \right)\ }+\mu \ {\sin  \left(\theta \right)\ }$ where 
$\theta $=\ artan $\frac{\left|Vec\left(q\right)\right|_Q}{Sc\left(q\right)}$,\ 0$\le \theta \le \pi $ and \ $\mu $ :=\ $\frac{Vec\left(q\right)}{\left|{  Vec}\left({  q}\right)\right|_Q}$ verifying ${\mu }^2 =\ -1$.

In this paper, we will study the  quaternion-valued signal $f:{\mathbb{R}}^2\to \mathbb{H} $ , $f$ which can be expressed as \\ \hspace*{4cm}$f=f_0+i f_1+jf_2+kf_3,$ \hfill(2.4) \\ with $f_m~: {\mathbb R^2}\to \ {\mathbb R}\ for\ m=0,1,2,3.$

Let the inner product of  ~$f,g: {\mathbb{R}}^2\to \mathbb{H} $ be defined by

$<f\left(x\right),g\left(x\right)>:=\int_{{\mathbb R}^2}{f}$(x)$\ \overline{g(x)}dx,$

we define ${|f|}^2_{2,Q}:=~<f,f>$\ and ${|f|}^2_{1,Q}= \int_{{\mathbb R}^2}{{|f(x)|}_Q}\ dx.$ \hfill(2.5) \\
Also we denote by ${L}^p \left({{\mathbb  R}}^{{  2}},{\mathbb  H}\right), p=1,2$ \ the space of integrable functions $f$ taking values in $\mathbb H$ such that ${\left|{  f}\right|}_{p,Q} <\infty $.

\section{The two-sided quaternion Fourier transform}
Ell[6] defined the quaternion Fourier transform (QFT) which has an important role in the representations of signals due to transforming a real 2D signal into a quaternion-valued frequency domain signal,  QFT belongs to the family of Clifford Fourier transformations because of (2.1).

There are three different types of QFT, the left-sided QFT , the right -sided QFT , and two-sided QFT [12].

We now review the definition and some properties of the two-sided QFT.
\begin{Def}
Let$\ f$ in    $L^1\left({\mathbb R}^2,{\mathbb H}\right)$. Then two-sided quarternionic Fourier transform of the function  $f$ is given by

${\mathcal F}\{f(x)\}$($\xi $)=$\int_{{\mathbb R}^2}{e^{-i{2\pi \xi }_1x_1}} f$($x$)$ e^{-j{2\pi \xi }_2x_2}dx,\ dx=dx_1dx_2.$ \hfill(3.1)

Where ~$\xi,x\in {\mathbb R}^2$
\end{Def}
We define a new module of $\mathcal F\{f\} $ as follows :

\hspace*{4 cm}${\left\|\mathcal F\left\{f\right\}\right\|}_Q\ :=\ \sqrt{\sum^{m=3}_{m=0}{{\left|\mathcal F\left\{f_m\right\}\right|}^2_Q}}. $ \hfill(3.2)

Furthermore, we define a new $L^2$-norm of ${\mathcal F\{f\}}$ as follows :\\
\hspace*{4 cm}$ {\left\|\mathcal F\{f\}\right\|}_{2,Q}:=\sqrt{\int_{\mathbb R^2}{{\left\|\mathcal F\left\{f\right\}(y)\right\|}^2_Qdy}}. $ \hfill(3.3)

It is interesting to observe that ${\left\|{\mathcal  F}\{f\}\right\|}_{Q}$    is not equivalent to$\ {\left|{\mathcal  F}\{f\}\right|}_{Q}$  unless  $f$ is real valued.\\

From (2.4) and the above defnition, we have the following lemma:

  \begin{Lemma}

Let $\ f\in L^1\left(\mathbb R^2,\mathbb H\right)$  then

\hspace*{6cm}$ {\mathcal F}\{ifj\} =\ i{\mathcal  F} \{f\}j.$ \hfill(3.4)

\end{Lemma}

\begin{Lemma}

  ${\left\|\ \right\|}_Q$  is indeed a norm.

Proof.  Verification of the positivity and homogeneity properties is straightforward, we will provide  the triangle inequality~

Let $f=f_0+if_1+jf_2+kf_3,\ g=g_0+ig_1+jg_2+kg_3\in L^1\left(\mathbb R^2,\mathbb H\right)$,

  By  linearity of  $\mathcal F$ we have

${ \mathcal F}\left\{{ f}\}+{ \mathcal F}{ \{}{ g}\right\}{ =}{ \mathcal F}\left\{{ f}{ +g}\right\}.$

By definition we obtain ${\left\|{ \mathcal F}\left\{{ f}\right\}+{ \mathcal F}\left\{{ g}\right\}\right\|}^2_Q{ =}\sum^{m=3}_{m=0}{{\left|\mathcal F\{f_m+g_m\}\right|}^2_Q}$

 Then, by  successive equivalences, we have

${\left\|{ \mathcal F}\left\{{ f}\right\}+{ \mathcal F}{ \{}{ g}\}\right\|}^2_Q\le \ {\left\|{ \mathcal F}\left\{{ f}\right\}\right\|}^2_Q+{\left\|{ \mathcal F}{ \{}{ g}\}\right\|}^2_Q$+2${\left\|{ \mathcal F}\left\{{ f}\right\}\right\|}_Q{\left\|{ \mathcal F}\left\{{ g}\right\}\right\|}_Q$

${{{\Leftrightarrow}}}\sum^{m=3}_{m=0}{{\left|\mathcal F\{f_m\}+{\mathcal F\{g}_m\}\right|}^2_Q}\le \sum^{m=3}_{m=0}{{\ (\left|\mathcal F\{f_m\}\right|}^2_Q}+{\left|\mathcal F\{g_m\}\right|}^2_Q)+$\textit{2}${\left\|\mathcal F\left\{f\right\}\right\|}_Q{\left\|\mathcal F\left\{g\right\}\right\|}_Q$\textit{)}

As ${\left|\mathcal F\left\{f_m\}+\mathcal F\{g_m\right\}\right|}^2_Q\le {\left|\mathcal F\left\{f_m\right\}\right|}^2_Q$\textit{+}${\left|\mathcal F\left\{g_m\right\}\right|}^2_Q$\textit{+2}${\ \left|\mathcal F\left\{f_m\right\}\mathcal F\left\{g_m\right\}\right|}_Q$

\hfill(triangle inequality for quaternion  norm)

$\sum^{m=3}_{m=0}{{\left|\mathcal F\left\{f_m\right\}\right|}_Q}{\left|\mathcal F\left\{g_m\right\}\right|}_Q\le {\left\|\mathcal F\left\{f\right\}\right\|}_Q{\left\|\mathcal F\left\{g\right\}\right\|}_Q$

${\Leftrightarrow}\sum^{m=3}_{m=0}{{\left|\mathcal F\left\{f_m\right\}\right|}_Q}{\left|\mathcal F\left\{g_m\right\}\right|}_Q\le \sqrt{\sum^{m=3}_{m=0}{{\left|\mathcal F\left\{f_m\right\}\right|}^2_Q}}\sqrt{\sum^{m=3}_{m=0}{{\left|\mathcal F\left\{g_m\right\}\right|}^2_Q}},$

which is a true proposition, according to the inequality of Cauchy-Schwarz inequality in $\mathbb R^4.$

\end{Lemma}

\begin{Lemma}{Inverse QFT} [3, Thm 2.5]

If  $f,{\mathcal F}\{f\}\in L\left({\mathbb{R}}^2,{\mathbb{H}}\right)$, then

$f(x)=\int_{{\mathbb{R}}^2}{e^{i2\pi {\xi }_1x_1}} {\mathcal F}\{f(x)\}(\xi ) e^{j{2\pi \xi }_2x_2}d^2\ \xi. $\hfill(3.5)
\end{Lemma}

\begin{Lemma}{Plancherel theorem for QFT} [4, Thm. 3.2]

If $f\in $ $L^2\left({\mathbb{R}}^2,{\mathbb H}\right)~$, then 
                               
${|f|}_{2,Q}={\left\|{\mathcal F}\left\{f\right\}\right\|}_{2,Q}.$\hfill(3.6)

\end{Lemma}

\begin{Lemma}
 ${\mathcal F}\left\{e^{-\pi {\left|x\right|}^2}\right\}\left(y\right)$=$e^{-\pi {\left|y\right|}^2}$
where ~$x,y \in {\mathbb R}^2.$\hfill(3.7)
 \end{Lemma}
Proof. ${\mathcal F}\left\{e^{-\pi {\left|x\right|}^2}\right\}\left(y\right)=\int_{{{ \mathbb R}}^2}{e^{-2\pi ix_1y_1 }e^{-\pi (x^2_1+x^2_2)}e^{-2\pi jx_2y_2}dx_1dx_2}$ 

\hspace*{3.9cm}=$e^{{ -}\pi {(y}^2_1+y^2_2)}\int_{{ \mathbb R}}{e^{-\pi  {(x_1+{i}y_1)}^2 }dx_1}\int_{{ \mathbb R}}{e^{-\pi  {(x_2+ { j}y_2)}^2 }dx_2}$ (Fubini)

We know that  $\int_{{ \mathbb R}}{e^{-{ z} {(t+z')}^2 }dt=}\sqrt{\frac{\pi }{z}}$,  for  ${ z},{ z'}\in {\mathbb C}$,\ Re(z)$>$0\ (Gaussian integral

with complex offset)\\ 
Therefore $\int_{{ \mathbb R}}{e^{-\pi  {(x_1+i y_1)}^2 }dx_1}=\int_{{ \mathbb R}}{e^{-\pi  {(x_2+ j y_2)}^2 }dx_2}$=1 which give us the desired result.

 \vspace*{0.5cm}Let\ *\ denote\ the\ convolution\ defined\ by\\
 For\ $ f,\ g \in \L^2 \left({\mathbb R}^2,{\mathbb  H}\right)\ \ {  f*g}\left(x\right):= \int_{{\mathbb R}^2}{{  f}\left({  t}\right){  g}\left({  x-t}\right)dt }.$
 \begin{Lemma}
 For $f,\ g\in  L^2\left({\mathbb R}^2,{\mathbb  R}\right) \ we\ have\ \  {\mathcal F}\{  f*g\}  ={\mathcal F}\{f \}\ {\mathcal  F}\{ g\}. $\hfill(3.8)
\end{Lemma}
Proof.\ See\ [1,\ Thm.\ 13].

\begin{Lemma}

Let $x\in {\mathbb R}$, then\\
$\frac{{\partial }^m}{{\partial x}^m}e^{-\pi x^2}=e^{-\pi x^2}P_m\left(x\right),  \ \ \  m \ge 0.  $  \hfill(3.9)

Where $P_m\left(x\right)$ is polynomial of degree m.
\end{Lemma}

Proof. I recall   that the polynomials $P_m$ are the  

Hermite functions modulo  $\pi ,$  ${(-1)}^m$;  defined for  $m\ge 0 $ by

\[H_m\left(x\right)   ={(-1)}^me^{-x^2}  \frac{{\partial }^m}{{\partial x}^m}e^{-x^2}.\] 

The proof  of lemma (3.8) is done by induction. It is trivial for $m = 0$.

Assume, for $m=k$, that (3.9) holds.

Let $m=k+1$. 

$\frac{{\partial }^{k+1}}{{\partial x}^{k+1}}e^{-\pi x^2}=$(-$2\ \pi xP_k\left(x\right)+P^{'}_k(x)$)$\ e^{-\pi x^2}$

\hspace*{1.8 cm}$ =e^{-\pi x^2}Q_{k+1}\left(x\right).$

\begin{Lemma}

  Let $\ f,\ x^m_1x^n_2f\in L\left({\mathbb R}^2,{\mathcal H}\right)$  for $(m,n)\ \in {\mathbb N}^2$ 

         Then

    ${\mathcal F}\{x^m_1x^n_2f\}(\xi )= {(\frac{1}{2\pi })}^{m+n} i^m\frac{{\partial }^{m+n}}{{\partial {\xi }_1}^m{\partial {\xi }_2}^n}{\mathcal F}\{f\}{(}\xi {)}j^n.$\hfill(3.10)
\end{Lemma}

Proof.\ See\ [5,\ Thm.\ 2.2 p.6].

\begin{Lemma}

Let  $f:{\mathbb{R}}^2\to \mathbb{R} $ be of the form

$f\left(x\right)=P(x)e^{-\pi {\left|x\right|}^2}$ where  $P$ is a polynomial

We have      ${\mathcal F}\{f\}(\xi )=\ Q(\xi )e^{-\pi {\left|\xi \right|}^2}$, 
where $Q$ is quaternion polynomial with $deg P =deg Q$.

\end{Lemma}

Proof. In two dimensions a complete l-th degree polynomial is given by  \\
  ${{P}}_{l}({{x}}_{{1}},{{x}}_{{2}})=\ \sum^{k=l}_{k=0}{{\alpha }_kx^m_1x^n_2} \  \ \   m+n \le  k$  .

We have by lemma (3.9)\ and (3.7)

${\mathcal F}\{x^m_1x^n_2e^{-\pi {\left|x\right|}^2}\} (\xi {)}=\ {(\frac{1}{2\pi })}^{m+n}\ i^m\frac{{\partial }^{m+n}}{{\partial {\xi }_1}^m{\partial {\xi }_2}^n}e^{-\pi {\left|\xi \right|}^2}j^n$

\hspace*{3.2 cm}=${(\frac{1}{2\pi })}^{m+n}\ i^mj^n\frac{{\partial }^m}{{\partial {\xi }_1}^m}e^{-\pi{\xi }^2_1}\frac{{\partial }^n}{{\partial {\xi }_2}^n}e^{-\pi{\xi }^2_2}$

\hspace*{3.2 cm}=${(\frac{1}{2\pi })}^{m+n}\ i^mj^nP_m\left({\xi }_1\right)Q_n\left({\xi }_2\right)e^{-\pi {\left|\xi \right|}^2},$  (by (3.9)).

The linearity of the two-sided QFT completes the proof.

\begin{Lemma}

Let $f:{\mathbb{R}}^2\to \mathbb{R} $  be of the form

$f\left(x\right)=P(x)e^{-\pi \alpha {\left|x\right|}^2}$ where  $P$ is a polynomial  and $\alpha >0$\\

Then      \[{\mathcal F}\{f\}(\xi {\rm )}=\ Q(\xi )e^{-\frac{\pi }{\alpha }{\left|\xi \right|}^2}\],

with $Q$ is quaternion polynomial with deg P =deg Q.
\end{Lemma}

Proof.  The proof is obtained by combining dilation property 

${\mathcal F}\{a f\}(\xi )={(\frac{1}{a})}^2\ {\mathcal F}\{ f(\frac{1}{a}\xi )\}, \ \   (a>0).\ \  $\ (See [3, Th 2.12]) , and lemma 3.10.

\section{BEURLING'S  THEOREM  FOR THE TWO-SIDED QUATERNION FOURIER TRANSFORM}
In this section, we provide Beurling's theorem for the two-sided quaternion Fourier transform. 
\begin{Lemma}
Let $f \in L^2(\mathbb R^2,\mathbb R) $\ and $d\ge 0$\ such that\\
\hspace*{4cm}$\int_{\mathbb R^2}{\int_{\mathbb R^2}{ \frac{{\left|{{  f}}\right|} \ {|\mathcal F\left\{{{  f}}\right\}\left(y\right)|}_Q}{{(1+\left|x\right|+\left|y\right|)}^d}{\ e}^{2\pi \left|x\right|\left|y\right|}}}\ dxdy <\infty, $\hfill(4.1)

 then  $f \in L^1\left({\mathbb R}^2,{\mathbb R}\right)$ and $\mathcal F\{f\}\in   L^1({\mathbb R}^2,{\mathbb H}).$

\end{Lemma}
Proof. We may assume that ${  f}\ne 0 $. By (4.1) and Fubini theorem, we obtain for almost every $y  \in \mathbb R^2$

\[{\left|\mathcal F\left\{{{  f}}\right\}\left(y\right)\right|}_Q\int_{\mathbb R^2}{\frac{\left|{  f}\right|   \ }{{(1+\left|x\right|+\left|y\right|)}^d}}e^{2\pi \left|x\right|\left|y\right|}\ dx <\infty. \] 

Since  ${  f}\ne 0$     we have\ by  the injectivity of\ ${  \mathcal F}$\ (Lemma 3.6),  \ ${  \mathcal F}\left\{{  f}\right\}\ne {  0}$, therefore there exists ${y}_0{  \ }\in $ ${{  \mathbb R}}^{{  2}}$, ${y}_0\ne {  0}$ such that ${  \mathcal F}\left\{{  f(}{y}_0{  )}\right\}\ne {  0}$, and

Therefore   $\int_{\mathbb R^2}{\frac{\left|{  f}\right|   \ }{{(1+\left|x\right|+\left|y_0\right|)}^d}}e^{2\pi \left|x\right|\left|y_0\right|}\ dx$ $<\infty. $ 

Since  $\frac{{{  e}}^{\left|x\right|\left|{y}_0\right|}   }{{(1+\left|x\right|)}^d}$  ${{  e}}^{\left|x\right|\left|{y}_0\right|}\ge {  1}$  for large $\left|x\right|$, it follows that $\int_{{{  \mathbb R}}^{{  2}}}{\left|{  f}\right|   \ }{  dx}$ ${  <}\infty $, so ${  f}$ $\in $ ${{  L}}^{{  1}}\left({{  \mathbb R}}^{{  2}},{  \mathbb R}\right)$.
 Interchanging the roles of $f$ and $\mathcal F\{f\}$, we get $ \mathcal F\{f\}\in {{  L}}^{{  1}}\left({{  \mathbb R}}^{{  2}},{\mathbb  H}\right)$ .

\begin{Th}

Let $f \in L^2\left(\mathbb R^2,\mathbb H\right){  \ }{  and\ \ d}\ge 0~$   satisfy  

\hspace*{4 cm} $\int_{\mathbb R^2}{\int_{\mathbb R^2}{\ \frac{{\left|{  f}\right|}_Q\  {\left\|\mathcal F\left\{f\right\}\left(y\right)\right\|}_Q}{{(1+\left|x\right|+\left|y\right|)}^d}{\ e}^{2\pi \left|x\right|\left|y\right|}}}\ dxdy$ $<\infty, $    \hfill(4.2)

Then $f(x)= P(x)\ e^{-a{\left|x\right|}^2}$ .

Where $a>0$ and $P$ is a polynomial of degree $<\frac{d-2}{2}$. In particular , $f$ is identically 0  when $d\le 2.$

\end{Th}

\begin{Rem}
It is important to see that for every component function ${{  f}}_{{  m}},m = 0, 1, 2, 3 $\ of the quaternion function $f$, we have by (2.2), (2.4) and (3.2)
 
 ${\left|{f}_{m}\right|} \le {\left|{  f}\right|}_Q$ and   ${\left|{{  {\mathcal  F}}{  \{}{  f}}_{m}\}\right|}_Q\le {\left\|{\mathcal  F}{  \{}{  f}\}\right\|}_Q$

So (4.2) implies $\int_{\mathbb R^2}{\int_{\mathbb R^2}{ \frac{{\left|{{  f}}_{m}\right|} \ {|\mathcal F\left\{{{  f}}_{m}\right\}\left(y\right)|}_Q}{{(1+\left|x\right|+\left|y\right|)}^d}{\ e}^{2\pi \left|x\right|\left|y\right|}}}\ dxdy <\infty $.\\
Then if the theorem is proved for ${{  f}}_{{  m}}\in L^2\left({\mathbb R}^2,{\mathbb R}\right)$ we obtain by (2.4) the result for \ $f$.\\

\end{Rem}

Proof. \textbf{First step.}\ Let $f \in L^2\left({\mathbb R}^2,{\mathbb R}\right)\ $ and suppose (4.1)

We define  $g=f*\varphi $  where $\varphi \left(x\right)=e^{-\pi {\left|x\right|}^2},x\in \mathbb R^2$.

It follows from (2.3), (3.7) and  (3.8) that\\
\hspace*{4cm}  ${\left|\mathcal F\left\{g\right\}(y)\right|}_Q ={\left|\mathcal F\left\{f \right\}\left(y\right)\right|}_Qe^{-\pi {\left|y\right|}^2}.$\hfill(4.3)\\
Then ${ \left|\mathcal F\left\{{ g}\right\}\left(y\right)\right|}_Q\le {\ \left|f\right|}_{1,Q}e^{-\pi {\left|y\right|}^2}$

and  by lemma 4.1 we obtain  $\mathcal F\left\{g\right\}\ e^{\pi {\left|.\right|}^2}\in L^1\left({ \mathbb R}^{{2}},{ \mathbb H}\right)$.\hfill(4.4)

We will show the following assumptions :

\hspace*{4 cm}$\int_{{  \mathbb R}^{{   2}}}{\int_{{  \mathbb R}^2}{\frac{\left|g(x)\right|\ {\left|{\mathcal   F}\left\{{   g}\right\}\left({   y}\right)\right|}_{{   Q}}\ }{{(1+\left|x\right|+\left|y\right|)}^d}}e^{2\pi \left|{   x}\right|\left|{   y}\right|}{   \ }dxdy} \ < \infty.$       \hfill(4.5)

For every $R >0$, there exists $C > 0$ such that 

\hspace*{3 cm}$\int_{\left|x \right|\le   \mathbb R}{\int_{{  \mathbb R}^{2}}{\left|{   g(x)}\right|{\ \left|{\mathcal    F}\left\{g\right\}\left(y\right)\right|}_{{   Q}}}{{   e}}^{{   2}\pi \left|{   x}\right|\left|y \right|}dy\ dx}\le \ {C\ (1+  R)}^d.$  \hfill  (4.6)

In view of (4.3) and the definition of $g$, the integral in (4.5)  is less than

\[\int_{{{  \mathbb R}}^{{  2}}}{\int_{{{  \mathbb R}}^{{  2}}}{\left|{  f(t)}\right|{\left|{  \mathcal F}\left\{{  f}\right\}\left(y\right)\right|}_{{  Q}}{  h(t,y)}}{{  e}}^{{  2}\pi \left|{  t}\right|\left|y\right|}{  \ dt\ dy}},\] 

where ${  h}\left({  t,y}\right){  \ }$:=${  \ }{{  e}}^{-\pi {\left|y\right|}^{{  2}}}{{  e}}^{-{  2}\pi \left|{  t}\right|\left|y\right|}\int_{{{  \mathbb R}}^{{  2}}}{\frac{{  1\ \ }}{{{  (1+}\left|x\right|+\left|y\right|{  )}}^{{  d}}}}{{{  e}}^{-\pi {\left|{  x-t}\right|}^{{  2}}}{  e}}^{{  2}\pi \left|x\right|\left|y\right|}{  \ dx}.$

To prove (4.5) we sould prove that ${  h}\left({  t,y}\right)\le {{  C\ (1+}\left|{  t}\right|+\left|y\right|{  )}}^{{-d}}.$

As $\left|{  <}x,t>\right|\le \left|x\right|\left|{  t}\right|$\ (Schwarz's inequality)  we have  
\[{  h}\left({  t}{  ,y}\right)\le \int_{{{  \mathbb R}}^{{  2}}}{\frac{{  1\ \ }}{{{  (1+}\left|x\right|+\left|y\right|)}^{d}}}{{  e}}^{-\pi {(\left|x\right|-(\left|{  t}\right|+\left|y\right|{  ))}}^2}{  \ dx}.\] 

Let $0<\delta <1$  and write A= ${  1+}\left|{  t}\right|+\left|y\right|$\\
 then\\
 $h\left({  t,y}\right)\le \int_{\left| \left|x\right|-\left(\left|{  t}\right|+\left|y\right| \right)\right| > \delta A\ }{\frac{{  1 }}{{(1+\left|x\right|+\left|y\right|)}^{d}}}{{  e}}^{-\pi {{  (}\left|x\right|-{  (}\left|{  t}\right|+\left|y\right|{  ))}}^2}{ dx}+\int_{\left|  \left|x\right|-\left(\left|{  t}\right|+\left|y\right|\right)\right| \le \delta A\ }{\frac{1}{{{  (1+}\left|x\right|+\left|y\right|{  )}}^{{  d}}}}{{  e}}^{-\pi {{  (}\left|x\right|-{  (}\left|{  t}\right|+\left|y\right|{  ))}}^2}dx$

It's clear that the first integral satisfies the  desired estimate, for the second integral,

by   the triangular inequality we have

$\left|{t}\right|{=}\left|\left|{t}\right|+\left(\left|{y}\right|{-}\left|{x}\right|\right){-}\left(\left|{y}\right|{-}\left|{x}\right|\right)\right|\le \left|\left|{t}\right|{+}\left|{y}\right|{-}\left|{x}\right|\right|{+}\left|\left|{y}\right|{-}\left|{x}\right|\right|$ \\
\hspace*{6.4cm}$\le \left|\left|{x}\right|-\left|{t}\right|{-}\left|{y}\right|\right|{+}\left|{y}\right|{+}\left|{x}\right|$

\hspace*{6 cm}$\le \left|\left|{x}\right|-\left|{t}\right|{-}\left|{y}\right|\right|{+}\left|{y}\right|{+2}\left|{x}\right|{+1},$ 

then \hspace*{2.8cm}${1}{+}{2}\left|{x}\right|{+}{2}\left|{y}\right|\ge  \left|{t}\right|{+}\left|{y}\right|{-}\left|\left|{x}\right|{-}\left|{t}\right|{-}\left|{y}\right|\right|,$ 

therefore    ${1+}\left|{y}\right|{+}\left|{x}\right|\ge \frac{{1}}{{2}}$+$\frac{\left|{t}\right|}{{2}}$+$\frac{\left|{y}\right|}{{2}}$-$\frac{\left|\left|{x}\right|{-}{(}\left|{t}\right|{+}\left|{y}\right|{)}\right|}{{2}}$

\hspace*{3.8 cm}$\ge \frac{1-\delta }{2} A, $  \hfill(because $\left| \ \left|x\right|-\left(\left|{  t}\right|+\left|y\right| \ \right| \right)\le \delta A$)

consequently the desired estimate is also satisfied by the second integral.

The proof of (4.5) is completed.

Now, we will demonstrate  (4.6),  choose  $\delta >1$  ,  we consider

$\int_{\left|{ x}\right|\le { R}}{\left|{ g}\left({ x}\right)\right|{ (}\int_{\left|{ y}\right|>2\delta R}{{\left|{ \mathcal F}\left\{{ g}\right\}\left({ y}\right)\right|}_{{ Q}}}{{ e}}^{2\pi \left|{ x}\right|\left|{ y}\right|}{ \ dy}+\int_{\left|{ y}\right|\le 2\delta R}{{\left|{ \mathcal F}\left\{{ g}\right\}\left({ y}\right)\right|}_{{ Q}}}{{ e}}^{2\pi \left|{ x}\right|\left|{ y}\right|}{ \ dy)dx}}.$ \hfill(4.7)

by combining (4.3)\ and lemma 4.1 we obtain \\
 \hspace*{4cm}${\left|{ \mathcal F}\left\{{ g}\right\}\left({ y}\right)\right|}_{{ Q}}\le $ C${{ e}}^{-\pi {\left|{ y}\right|}^2}.$\hfill(4.8)\\
If $\left|{ x}\right|\le { R}<\frac{{ 1}}{2\delta }\left|{ y}\right|$ we have $2\pi \left|{ x}\right|\left|{ y}\right|$ $\le \frac{\pi }{\delta }{\left|{ y}\right|}^2$.\\
  As a consequence of $\int_{\left|{ y}\right|>2\delta R}{}{{ e}}^{-\pi { (1-}\frac{{ 1}}{\delta }{ )}{\left|{ y}\right|}^2}{ \ dy\ }<\infty $, we have    
\[\int_{\left|{ x}\right|\le { R}}{\left|{ g}\left({ x}\right)\right|{ (}\int_{\left|{ y}\right|>2\delta R}{{\left|{ \mathcal F}\left\{{ g}\right\}\left({ y}\right)\right|}_{{ Q}}}{{ e}}^{2\pi \left|{ x}\right|\left|{ y}\right|}{dy)dx}} \le { C}{\left|{ g}\right|}_{{ 1,Q}}.\]
On the other hand, if we multiply and divide by $(1+|x|+|y| )^d$ in the integral of right side in (4.7), we get\\
 $\int_{\left|x\right|\le R}{\int_{\left|y\right|\le 2\delta R}{\left|{ g}\left(x\right)\right|{\left|{ \mathcal F}\left\{{ g}\right\}\left({ y}\right)\right|}_{{ Q}}}{{ e}}^{2\pi \left|{ x}\right|\left|{ y}\right|}{ \ dy}{ dx}}\le {\left({ 1+}{ R}\right)}^d \int_{{\mathbb R}^2}{\int_{{\mathbb R}^2}{\frac{\left|{ g}\left({ x}\right)\right|{\left|{ \mathcal F}\left\{{ g}\right\}\left({ y}\right)\right|}_{{ Q}}}{{\left({ 1+}\left|{ x}\right|+\left|{ y}\right|\right)}^d}} e^{2\pi \left|{ x}\right|\left|{ y}\right|}dxdy}$\\
\hspace*{7.9 cm}$\le {C \left(1+R\right)}^d.$

This proves (4.6).

\textbf{Second step.} By lemma 3.4 we have\ ${{ g(x)=\mathcal F}\left\{{  \mathcal F}\left\{g\right\}\right\}\left(-{  x}  
  \right) \forall x \in {\mathbb R}^2}$, 
Furthermore 

Complexifying the variable$\ \ z=a+i_{\mathbb C}{  \ b\ };{  \ a}=(a_1,a_2$),${  \ b}=(b_1,b_2)\ \in {\mathbb R}^2$ (we note by $i_{\mathbb C}$ the complex number checking $i^2_{\mathbb C}$= -1)

We have 

\[g \left(z\right)=\int_{{{\mathbb R}}^2}{e^{2\pi i y_1(a_1+i_{\mathbb C} b_1)}{\mathcal F}\left\{g\right\}\left(y\right) e^{2\pi j y_2(a_2+i_{\mathbb C}b_2)}dy_1dy_2}\] 

then 
\[{|g(z)|}_Q\le \int_{{{\mathbb R}}^2}{{|{\mathcal F}\left\{g\right\}\left(y\right)|}_Q}e^{2\pi {  (}\left|y_1a_1\right|+\left|y_1b_1\right|{  +}\left|y_2a_2\right|+\left|y_2b_2\right|}dy_1dy_2\] 
By (4.8) \ \ \ \ ${|g\left(z\right)|}_Q\le Ce^{\pi {  (}{\left|a_1\right|^2+\left|b_1\right|^2)}}e^{\pi {  (}{\left|a_2\right|^2+\left|b_2\right|^2)}}\int_{{{\mathbb R}}^2}{e^{-\pi {  (}{\left|y_1\right|-(\left|a_1\right|+\left|b_1\right|))}^2}}e^{-\pi {  (}{\left|y_2\right|-(\left|a_2\right|+\left|b_2\right|))}^2}dy_1dy_2$

\hspace*{2.7cm}= $Ce^{\pi {\left|z\right|}^2}\int_{{\mathbb R}}{e^{-\pi {  (}{\left|y_1\right|-(\left|a_1\right|+\left|b_1\right|))}^2}}dy_1\int_{{\mathbb R}}{e^{-\pi {  (}{\left|y_2\right|-(\left|a_2\right|+\left|b_2\right|))}^2}}dy_2$\\
Since $\int^{+\infty }_{-\infty }{e^{-\pi {  (}{\left|t\right|+m)}^2}}\ dt=\int^{+\infty }_0{e^{-\pi {  (}{t+m)}^2}}\ dt$ +$\int^0_{-\infty }{e^{-\pi {  (}{-t+m)}^2}}\ dt\ $ for $m\in {\mathbb R}$

\hspace*{4cm}=$\int^{+\infty }_0{e^{-\pi {  (}{t+m)}^2}}\ dt$+$\int^{+\infty }_0{e^{-\pi ({t+m)}^2}}\ dt$\\
\hspace*{4.4cm}=2$\int^{+\infty }_m{e^{-\pi t^2}} dt$\\
\hspace*{4.4cm}$ \le 2 \int^{+\infty }_{-\infty }{e^{-\pi t^2}}\ dt=\ 2.$

We deduce that ${|g\left(z\right)|}_Q\le 4 \ Ce^{\pi {\left|z\right|}^2}.$\\
It follows that g is entire of order 2.\\
  
  \textbf{Third step.} The function $g$ admits an holomorphic extension to ${\mathbb C}^2$ that is of order 2. Moreover, there exists a polynomial R such that for all $z \in {\mathbb C}^2$,

$g(z)g(i_{\mathbb C}z) = R(z)$.

For all $x \in{\mathbb R}^2$ and $\theta \in {\mathbb R}$ , ${\left|{ g(}{{ e}}^{{ i}_{\mathbb C}\theta }{ x)}\right|}_Q \le  \int_{{\mathbb R}^2}{{ \left|{\mathcal F}\left\{g\right\}\left({ y}\right)\right|}_{{ Q}}}{{ e}}^{{ 2}\pi \left|{ x}\right|\left|{ y}\right|}{ \ dy}.$

We should prove that    ${ g(z)g(i_{\mathbb C}z)}  for\ z \in{\mathbb C}^2{ \ }$,   is a polynomial

To show that, we  define a new function $G$ on ${\mathbb C}^2$ by :\ $G : z \to \int^{z_1}_0{\int^{z_2}_0{{ g(u)g(iu)du}}}$.

$G$ is entire of order 2, because g it is. 

As $g\left(z\right){ g}\left({ i_{\mathbb C}z}\right)=\frac{\partial }{\partial z_1}\frac{\partial }{\partial z_2}$ G($z$), to prove our claim it is enough to show that $G$ is a polynomial, For this we use(4.6), and  the  Phragm\`en-Lindelh\"of’s principle[8], by folollowing the proof of [2, Prop. 2.2  page 32-33] we find that $g(y)=P(y){{ e}}^{({ By,y}{ )}}$ where $B$ a is symmetric matrix and $P(y)$ is a polynomial.
A direct computation shows that the form of the matrix   $B$ imposed by condition (4.5) is   $B=\delta I.$

Then $g\left(y\right)=$ $P(y)e^{-\delta {\left|y\right|}^2}$ , therefore by lemma (3.11)\ ${\mathcal F}\{g\}$  has a similar form.

As ${\mathcal F}\left\{{g}\right\}\left(y\right)={\mathcal F}\left\{{f}\right\}\left(y\right)e^{-\pi {\left|y\right|}^2}$, $f$ will be as in the theorem.$ \blacksquare$

\section{APPLICATIONS TO OTHER UNCERTAINTY PRINCIPLES}
In this section we derive some other versions of uncertainty principle for the two-sided quaternion Fourier transform.

\begin{Cor}
\textbf{(Hardy type)}\\
Let 
 ${ f}\in {{ L}}^{2}\left({\mathbb R}^{2}{ , {\mathbb H}}\right) \ and\  d\ge { 0}, \  \alpha { ,}\beta >0$.  with 

\[{\left|{ f}\right|}_Q\le C{{ (1+}\left|{ x}\right|{ )}}^{{ d}}{e}^{-\pi \alpha {\left|x\right|}^2}, \ \ {\left\| \ {\mathcal F} \left\{f\right\}\left(y\right)\right\|}_Q\le C{{ (1+}\left|{ y}\right|{ )}}^{{ d}}e^{-\pi \beta {\left|y\right|}^2}.\] 

(i)  If $\alpha \beta > \ 1 $, then $f=0 $.\\
\ (ii) If $\alpha \beta ={ 1}$, then $\ f(x) =P(x)\ {e}^{{ -}\pi \alpha {\left|x\right|}^2}$ , where  $P$ is a polynomial of degree $\le d$.

(iii) else there are infitely many linearly independent functions satisfying the conditions.

\end{Cor}

Proof. Firstly, from the remark (4.3)  it is enough to show the corollary  for $f \ \in {{ L}}^{{ 2}}\left({{  \mathbb R}}^{{ 2}},{  \mathbb R}\right)$.

 Form the decay conditions we have $\int_{ \mathbb R^2}{\int_{ \mathbb R^2}{\ \frac{{\left|{ f}\right|}_Q \ {\left\|{\mathcal F}\left\{f\right\}\left(y\right)\right\|}_Q}{{(1+\left|x\right|+\left|y\right|)}^d}{\ e}^{2\pi \left|x\right|\left|y\right|}}}\ dxdy,$

is bounded  by a constant multiple of 

$\int_{ \mathbb R^2}{\int_{ \mathbb R^2}{\ \frac{{{ (1+}\left|{ x}\right|{ )}}^{d}\  {{ (1+}\left|{ y}\right|{ )}}^{{ d}}}{{(1+\left|x\right|+\left|y\right|)}^d}{\ e}^{-\pi [\alpha { \ }{\left|x\right|}^2+\beta { \ }{\left|y\right|}^2-2\left|x\right|\left|y\right|]}}}\ dxdy.$  \hfill(5.1)

Hence if $\beta >\frac{1}{\alpha }$,    (5.1)   is finite, because $\alpha { \ }{\left|x\right|}^2+\beta { \ }{\left|y\right|}^2>{ \ }\alpha { \ }{\left|x\right|}^2+\frac{1}{\alpha }{\left|y\right|}^2\ge 2\left|x\right|\left|y\right|.$

 So by Theorem (4.2)\ $f\left(x\right)= P\left({ x}\right){{ e}}^{{ -}\pi \alpha {\left|x\right|}^2}.$

From the decay condition on $ f$ it is clear that $ deg P \le d$.

But then ${\mathcal F}\left\{f\right\}\left(y\right)=Q(y){ \ }{{ e}}^{{ -}\frac{\pi }{\alpha }{\left|y\right|}^2},$   ( by lemma  3.11)

which cannot satisfy the decay condition  in the corollary as 

$\beta -\frac{1}{\alpha }>0$, then  $f=0$.

If $\beta =\frac{1}{\alpha }$, then  $f\left(x\right){ =P}\left({ x}\right){{ e}}^{{ -}\pi \alpha {\left|x\right|}^2}$with $deg P \le d.$

For the case $\alpha \beta <1$, let $\delta $ be such that $\beta <\frac{1}{\delta }<\frac{1}{\alpha }$.\ The function

$g\left(x\right){ =}{  R}\left({ x}\right){{ e}}^{-\pi \delta {\left|x\right|}^2}$\ will satisfy the conditions of the corollary with ${ R}\left({ x}\right){ \ }$is any polynom of degree less than $d$.

\begin{Rem}

The above corollary  is a generalization of the theorem [10, Thm 5.3].
\end{Rem}

\begin{Cor}
\textbf{(Gelfand-Shilov type)}

Let ${ \ and \ d}\in {\mathbb N},\ ~$   $\alpha { ,}\beta \ge 1$, ${ 1}<p,q< \ \infty \ with\ 1/{{ p}}{ +}{{ 1}}/{{ q}}{ =1\ \ }.$

Assume that ${ f}$ $\in $ ${{ L}}^{{ 2}}\left({\mathbb R}^{2},{\mathbb H}\right)$ satisfies

$\int_{{\mathbb R}^2}{\ \frac{{\left|{ f}\right|}_{{ Q}}{ \ \ }}{{{ (1+}\left|{ x}\right|{ )}}^{{ d}}}{{ \ e}}^{{ 2}\pi \frac{{\alpha }^{{ p}}}{{ p}}{\left|{ x}\right|}^{{ p}}}}{ \ dx} < \ \infty ,\ \ \int_{{\mathbb R}^{{ 2}}}{\frac{{ \ \ }{\left\|{\mathcal F}\left\{{ f}\right\}\left({ y}\right)\right\|}_{{ Q}}}{{(1+\left|{ y}\right|{ )}}^{{ d}}}{{ \ e}}^{{ 2}\pi \frac{{\beta }^{{ q}}}{{ q}}{\left|y\right|}^{{ q}}}}{ \ dy} < \ \infty. $

Then

(i)$ f=0  \ if \ (p,q) \ne \left(2,2\right)$   or   $\alpha \beta >1$.

(ii)   else , $f(x)= P(x){{ e}}^{{ -}\pi{\alpha }^{{ 2}}{\left|{ x}\right|}^{{ 2}}}$, where  $P$ is a polynomial of degree $< d-2$.

\end{Cor}

Proof. From the well-known Young's inequality $\xi \lambda { \ } \le ({\xi }^p /p) + ({\lambda }^q /q)$

 , valid for nonnegative real numbers $\xi $ and $\lambda $, we have that

$\alpha \beta \textbar x\textbar \textbar y  \le ({\alpha }^p/p){{ |}x{ |}}^p + ({\beta }^q/q){{ |}y{ |}}^q,$

hence the integral

\[\int_{\mathbb R^2}{\int_{\mathbb R^2}{\ \frac{{\left|{ f}\right|}_Q \ {\left\|{\mathcal F}\left\{f\right\}\left(y\right)\right\|}_Q}{{{ (1+}\left|{ x}\right|{ )}}^{{ d}}{{ (1+}\left|{ y}\right|{ )}}^{{ d}}}{\ e}^{2\pi \alpha \beta \left|x\right|\left|y\right|}}}\ dxdy,\] 

is finite, because it is bounded by

\[\int_{\mathbb R^2}{\int_{\mathbb R^2}{ \ \frac{{\left|{ f}\right|}_Q {\left\|{\mathcal F}\left\{f\right\}\left(y\right)\right\|}_Q}{{{ (1+}\left|{ x}\right|{ )}}^{{ d}}{{ (1+}\left|{ y}\right|{ )}}^{{ d}}}{\ e}^{{ 2}\pi \frac{{\alpha }^p}{p}{{ |}x{ |}}^p+{ 2}\pi { \ }\frac{{\beta }^q}{q}{{ |}y{ |}}^q}}}\ dxdy.\] 

So by Theorem (4.2)  it follows that $f=0$ \   if \  $\alpha \beta >1$.

And if  $\alpha \beta =1,\ f(x)= P(x){ \ }{{ e}}^{{ -}\gamma {\left|{ x}\right|}^{{ 2}}}$ where $\gamma >0\ $  with deg P $<$ d-2.

But if $p>2$,\ then $\int_{{{\mathbb  R}}^{{ 2}}}{\frac{{\left|{ f}\right|}_{{ Q}}{ \ \ }}{{{ (1+}\left|{ x}\right|{ )}}^{{ d}}}{{ \ e}}^{{ 2}\pi \frac{{\alpha }^{{ p}}}{{ p}}{\left|{ x}\right|}^{{ p}}}}{ \ dx}$ cannot be finite, similarly the decay condition on ${\mathcal F}\left\{f\right\}$ cannot be finite if\ $q>2$. So we must have $p=q=2$.

\begin{Cor}

\textbf{(Cowling-Price type)}
Let $f \in {L}^{2}\left({\mathbb R}^{2}{ , {\mathbb H}}\right) \ and \ d\ge { 0}~$   satisfy  

\[\int_{{\mathbb R}^{2}}{{ \ }{(\frac{{\left|{ f}{ (x)}\right|}_Q}{{{ (1+}\left|{ x}\right|)}^d}{e}^{2\pi {\alpha \left|{ x}\right|}^{2}}}}{ )}^p{ \ dx} \ { <}\ \infty ,\ \int_{{\mathbb R}^{2}}{{(\frac{{\left\|{\mathcal F}\left\{f\right\}\left(y\right)\right\|}_Q{ \ \ }}{{{ (1+}\left|{ y}\right|{ )}}^{{ d}}}{{ \ e}}^{2\pi {\beta \left|{ y}\right|}^{2}})^q}}\ dy \ { <}\ \infty, \]

with  ${ 1}<p,q<\infty ,\ {{ 1}}/{{ p}}{ +}{{ 1}}/{{ q}}{ =1.\ }$ \\Then 

(i)  $f=0,$ \ if   $\alpha \beta > \frac{1}{4}. $  

(ii)  $f(x)= P(x){ \ }e^{-2\pi\alpha {\left|x\right|}^{2}}$\ if $\alpha \beta =\frac{1}{4}, \ where \  $P$\  is \ a \ polynomial \ of \ degree \ $<$\  min\{(d - 2)/p,(d - 2)/q\} .$

\end{Cor}

Proof.\ By H\^older's inequality, we have $\int_{{{\mathbb R}}^{{ 2}}}{\int_{{\mathbb R}^2}{ \ }\frac{{\left|{ f(x)}\right|}_Q{ \ \ }}{{{ (1+}\left|{ x}\right|{ )}}^{{ d}}}}{ \ }\frac{{\left\|{\mathcal F}\left\{f\right\}\left(y\right)\right\|}_Q{ \ \ }}{{{ (1+}\left|{ y}\right|{ )}}^{{ d}}}{{ \ e}}^{{ 2}\pi {\alpha \left|{ x}\right|}^{{ 2}}{ +2}\pi {\beta \left|{ y}\right|}^{{ 2}}}{ dx}$ ${ \ \ dy}$ ${ <}\infty. $

So by taking $\alpha { =}\frac{{{ m}}^{{ 2}}}{2}$, $\beta { =}\frac{{{ n}}^{{ 2}}}{2}$  ($\alpha \beta { =}\frac{1\ }{4}\ when\ mn=1$)  we see that the result is a particular case of the previous corollary .

%\textbf{References}\\

\end{document}